\documentclass[12pt]{amsart}
\usepackage{amscd, amssymb,supertabular,graphicx}

\theoremstyle{plain}
\newtheorem{theorem}{Theorem}[section]
\newtheorem{lemma}[theorem]{Lemma}
\newtheorem{question}[theorem]{Question}
\theoremstyle{definition}

\newtheorem{remark}[theorem]{Remark}
\newtheorem{coro}[theorem]{Corollary}
\newtheorem{prop1}[theorem]{Proposition}
\newtheorem*{acknowledgment}{Acknowledgments}
\newtheorem*{convention}{Convention}

\theoremstyle{definition}

\begin{document}

\title{Set theoretic complete intersection for curves in a smooth affine algebra}

\begin{abstract}
Let $A$ be a regular ring of dimension $d$ ( $d \ge 3$ ) containing an infinite field $k$. Let $n$ be an integer such that $2n \ge d+3 $. Let $I$ be an ideal in $A$ of height $n$ and $P$ be a projective $A$-module of rank $n$. Suppose $P \oplus A \approx A^{n+1}$ and there is a surjection $\alpha$: $P \to I$. It is proved in this note that $I$ is a set theoretic complete intersection ideal. As a consequence, a smooth curve in a smooth affine $\mathbb{C}$-algebra with trivial conormal bundle is a set theoretic complete intersection if its corresponding class in the Grothendieck group is torsion. 
\end{abstract}

\date{June 8, 2005. Revised: August 10, 2005}

\curraddr{Department of Mathematics, Washington
University in St. Louis, One Brookings Dr. Campus Box 1146, St. Louis, MO 63130, USA}
\author{ZE MIN ZENG}
\email{zmzeng@math.wustl.edu}
\maketitle
\thispagestyle{empty}

\par

\section{Introduction}

Let $A$ be a commutative Noetherian ring of dimension $d$ ( $d \ge 3$ ). Let $J$ be a local complete intersection ideal in $A$ of height $d-1$. Then by the well known Ferrand-Szpiro construction[\ref {Sz} ], there exists a local complete intersection ideal $I$ which is contained in $J$, such that $\sqrt{I} = \sqrt{J}$ and $I/{I^2}$ is free $A/I$-module of rank $n$. One can ask when is such an ideal $I$ a set theoretic complete intersection? Inspired in part by the results in [ \ref {Bo}], [\ref {Lg} ], [\ref {Mu1}] and [\ref {Mu2}], it is conceivable that the property of an ideal being a set theoretic complete intersection is related to its class in the Chow group or Grothendieck group being a torsion element, thus it is natural to ask the following question:

\begin{question} \label{qu0}
Let $A$ be a commutative Noetherian ring of dimension $d$ ( $d \ge 3$ ). Let $I$ be a local complete intersection ideal in $A$ of height $n=d-1$ such that $I/{I^2}$ is a free $A/I$-module of rank $n$. Suppose $(A/I)$ is torsion in $K_0(A)$. Is $I$  a set theoretic complete intersection in $A$?
\end{question}

When $n \ge 5 $  and odd, the above Question has an affirmative answer (for example, see the proof of  [\ref {zzm2}, Theorem 3.2] ). For the case when $n$ is even, the author gave an affirmative answer to the above Question too in  [\ref {zzm2}, Theorem 3.6] if $A$ is a polynomial algebra containing $\mathbb{Q}$. This note is another attempt by the author to settle the above Question.

\par
In this note, we shall give an affirmative answer to Question \ref {qu0} in the case when $A$ is a smooth affine $\mathbb{C}$-algebra( see Corollary \ref {cor3}).

\begin{theorem} \label{thm0}
Let $A$ be a smooth affine $\mathbb{C}$-algebra of dimension $n+1 $, where $n \ge 4$ and even. Let $I$ be a local complete intersection ideal of height $n$ such that $I/I^2$ is a free $A/I$-module of rank $n$. Suppose $(A/I)$ is torsion in $K_0(A)$. Then $I$ is a set theoretic complete intersection in $A$.
\end{theorem}

All rings in this paper are assumed to be commutative and Noetherian. 
All modules considered are assumed to be finitely generated. We denote by $K_0(A)$ the Grothendieck group of projective modules over the ring $A$.

\section{Main theorem}

Let $A$ be a commutative Noetherian ring of dimension $n$, and $I \subseteq A$
be a local complete intersection of height $r$ ($ r \le n $). Suppose
 $I/I^2$ is $A/I$-free with base $\bar{f_1},\dots, \bar{f_r}$, $f_i
\in I$, $\bar {f_i}$ is the class of $f_i$ in $I/I^2$. Let 
$J=I^{(r-1)!} + (f_1, \dots, f_{r-1})$. Then, by a result of Murthy [\ref {Mu1}, Theorem 2.2],  there exists a surjection $P \to J $ with $P$ a projective $A$-module of rank $r$, such that 
$(P)- (A^r)= -(A/I) $ in $K_0(A)$.  Therefore to show Question \ref {qu0} has an affirmative answer in the case when $A$ is a smooth affine $\mathbb{C}$-algebra, it suffices to answer the following much more general question positively:
 
\begin{question} \label{qu}
Let $A$ be a regular ring of dimension $d$ ($d \ge 3$) containing an infinite field $k$. Let $n$ be an integer such that $2n \ge d+3 $. Let $I$ be an ideal in $A$ of height $n$ and $P$ be a projective $A$-module of rank $n$. Suppose $P \oplus A \approx A^{n+1}$ and there is a surjection $\alpha$: $P \to I$. Is $I$ a set theoretic complete intersection ideal in $A$ ?
\end{question}

\begin{remark} \label{rmk}
If $n$ is odd, we have the following proposition:
\par
\begin{prop1} \label{prop0}
Let $A$ be a ring of dimension $d$ ($d \ge 3$) and $n$ be a odd integer such that $2n \ge d+3 $. Let $I$ be a local complete intersection ideal in $A$ of height $n$ and $P$ be a projective $A$-module of rank $n$. Suppose $P \oplus A \approx A^{n+1}$ and there is a surjection $\alpha$: $P \to I$. Then $I$ is a set theoretic complete intersection ideal in $A$.
\end{prop1} 

\begin{proof}
Since $n$ is odd, $P$ has a free summand of rank 1 by Bass [\ref {Ba1}], say $P \approx Q \oplus A $. Let $x$ be the image of $(0, 1)$ under $\alpha$ and $J$ be the image of $Q$ under $\alpha$, then $I = (J, x)$. By some suitable elementary transformations on $P$, we may assume ht$J=n-1$. $2n \ge d+3$ implies rank($ Q/{JQ} $) $\ge$ dim$(A/J) +1 $ . By Bass cancellation [\ref {Ba2}], $Q/{JQ}$ is a free $A/J$-module of rank $n-1$. By [\ref {Mk}, lemma 1], $I$ is generated by $n$ elements. In particular,  $I$ is a set theoretic complete intersection ideal in $A$. The proof of the propositon is complete.
\end{proof}
Therefore, if $n$ is odd, Question \ref {qu} has an affirmative answer.
\end{remark}

In order to give a complete answer to the Question \ref {qu}, we need a few lemmas.

First, let us restate a lemma of Van der Kallen [ \ref {Vk1}, Lemma 4.9]:
\begin{lemma} \label{lm1}
Let $A$ be a commutative ring.  Let $[a_0, a_1,\dots , a_n]$ be a unimodular row over $A$ and $P$ be the cokernal of the natural map: 
 \[ \begin{CD} A@>(a_0^2, a_1, \dots, a_n)>>A^{n+1} . \end{CD} \]

Then the projective $A$-module $P$ has a free summand of rank 1.
\end{lemma}

Let $A$ be a regular ring of dimension $d$ ($d \ge 3$) containing an infinite field $k$. Let $n$ be an integer such that $2n \ge d+3 $. Let $[a_0, a_1,\dots , a_n] \in Um_{n+1}(A)$, and $P=A^{n+1}/{\sum_{0 \le i \le n } a_ie_i }$. The Euler Class group of $A$ is defined by Bhatwadekar and Sridharan in [\ref {BS1}], which is denoted by $E^n(A)$. And they also attached $P$ an element $e([a_0, a_1,\dots , a_n]) $ in $E^n(A)$.

The following is an easy corollary to the above lemma:

\begin{coro} \label{cor1}
Let $A$ be a regular ring of dimension $d$ ($d \ge 3$) containing an infinite field $k$. Let $n$ be an integer such that $2n \ge d+3 $. Suppose $[a_0, a_1,\dots , a_n] \in Um_{n+1}(A)$. Then $e([a_0^2, a_1,\dots , a_n]) = 0$ in $E^n(A)$.
\end{coro}

\begin{proof}
By [ \ref {BS1}, Theorem 5.4] and the above lemma, we are done.
\end{proof}

Let $B$ be a ring and $M$ be a finitely generated $B$-module. We use the following convention throughout the rest of this note for simplicity:

\begin{convention}
Let $m_1, \dots, m_n \in M$. We say a map: $B^n \to M $ is given by ( $m_1, \dots, m_n$ ) to mean a $B$-module homomorphism: $B^n \to M $ defined by sending $e_i$ to $m_i$ for $i=1, \dots, n$, where $(e_1, \dots, e_n)$ is a standard basis of $B^n$.
\end{convention}
 
\begin{lemma} \label{lm2}
Let $A$ be a regular ring of dimension $d$ ($d \ge 3$) containing an infinite field $k$ and $n$ be an even integer such that $2n \ge d+3 $. Let $I=(f_1, \dots, f_n) $ be an ideal of height $n$ and $u \in A $ such that $1-uv \in I$. Let  $ u^2\omega$  be the surjection: $(A/I)^n \to I/I^2$ given by $(\overline{u^2f_1}, \overline{f_2}, \dots, \overline{f_n} )$, where bar denotes reduction modulo $I$. Then $(I, u^2\omega)=0$ in $E^n(A)$.
\end{lemma}

\begin{proof}
Applying Lemma 5.6 in [\ref {BS1}], we have $e([v^2, f_1, \dots, f_n])= (I, u^2\omega)$ in $E^n(A)$. By Corollary \ref {cor1}, we are done.
\end{proof}

The proof of the following lemma is analogous to the proof of Lemma 5.4 in [\ref {BS2}]. 
\begin{lemma} \label{lm3}
Let $A$ be a regular ring of dimension $d$ ($d \ge 3$) containing an infinite field $k$ and $n$ be an even integer such that $2n \ge d+3 $. Let $I$ be an ideal of height $n$ and $u \in A $ such that $1-uv \in I$. Suppose $I=(f_1, f_2, \dots,f_n) + I^2$. Let $ \omega$  be the surjection: $(A/I)^n \to I/I^2$ given by $( \overline{f_1},\overline{f_2}, \dots, \overline{f_n} )$ and let $ u^2\omega$  be the surjection: $(A/I)^n \to I/I^2$ given by  $( \overline{u^2f_1},\overline{f_2}, \dots, \overline{f_n} )$ . Then $(I, \omega)=(I, u^2\omega)$ in $E^n(A)$.
\end{lemma}

\begin{proof}
If $(I, \omega)=0$ in $E^n(A)$, then we are done by Lemma \ref {lm2}. So we may assume $(I, \omega)\ne 0 $, then by corollary 2.4 in [3] we can find an ideal $J$ of height $n$, such that $I+J=A$, $J \cap I = (f_1, \dots, f_n)$ and $J=(f_1, \dots, f_n) + J^2$. Let $\omega_J$ be the surjection: $(A/J)^n \to J/J^2$ given by $( \overline{f_1},\overline{f_2}, \dots, \overline{f_n} )$, then $(I, \omega)+(J, \omega_J)=0 $ in $E^n(A)$. Since $I+J=A$, we can write $1-u=x+y$ for some $x \in I$, $y \in J$. Let $b=1-y$, then $b = 1$ modulo $J$ and $b = u$ modulo $I$. By Lemma \ref {lm2} above and Theorem 4.2 in [ \ref {BS1}], we see that there exists a surjection $ \phi $: $A^n \to I \cap J $, such that $\phi \otimes A/I$ is given by $( \overline{u^2f_1},\overline{f_2}, \dots, \overline{f_n} )$  and $ \phi \otimes A/J $ is given by $( \overline{f_1},\overline{f_2}, \dots, \overline{f_n} )$. From the surjection $ \phi $, we have $(I, u^2\omega)+(J, \omega_J)=0 $ in $E^n(A)$. Combining the relation $(I, \omega)+(J, \omega_J)=0 $, we have $(I, \omega)=(I, u^2\omega)$ in $E^n(A)$.
\end{proof}

Let $A$ be a regular ring of dimension $d$ ($d \ge 3$) containing an infinite field $k$ and $n$ be an integer such that $2n \ge d+3 $. By a theorem of Van der Kallen [ \ref {Vk2}, Theorem 4.1], the universal weak Mennicke symbol $$wms: Um_{n+1}(A)/{E_{n+1}(A)} \to WMS_n(A)$$ is a bijection with an abelian target, which provides $Um_{n+1}(A)/{E_{n+1}(A)}$ with the desired structure of an abelian group. In [\ref {BS1} , Theorem 5.7], Bhatwadekar and Sridharan showed that the natural map  $$ e: Um_{n+1}(A)/{E_{n+1}(A)} \to E^n(A) $$ is a group homomorphism, where the group structure of  $ Um_{n+1}(A)/{E_{n+1}(A)} $ is the one defined above by Van der Kallen. 

\par

Let $[a_0, a_1, \dots, a_n]$ be a unimodular row over $A$ and let $$P=A^{n+1}/{\sum_{0 \le i \le n } a_ie_i },$$ where $(e_0, \dots, e_n)$ is a standard basis of $A^{n+1}$. Let $p_i$ denote the image of $e_i$ in $P$, then $P= \sum_{0 \le i \le n } Ap_i $ and $\sum_{0 \le i \le n } a_ip_i=0$. Suppose there is a surjection $\alpha: P \to I $, where $I$ is an ideal of height $n$. Let $f_i$ be the image of $p_i$ under the surjection $\alpha$. Since $n+1 \ge $ dim $(A/I) + 2$, $[\overline{a_0},\overline{a_1}, \dots, \overline{a_n}]$ is completable to an elementary matrix in $E_{n+1}(A/I)$. So we may assume $[a_0, a_1, \dots, a_n] \equiv [1, 0,\dots, 0]$ modulo $I$. $2n \ge d+3 $ implies $P/IP$ is a free $A/I$-module of rank $n$ by Bass cancellation.  Since $\sum_{0 \le i \le n } a_ip_i=0$, we can write $ I= (f_1, \dots, f_n) + I^2 $ . So if we let $\omega : (A/I)^n \to I/{I^2} $ denote the surjection given by $( \overline{f_1}, \dots, \overline{f_{n-1}}, \overline{f_n} )$, then $e([a_0, a_1, \dots, a_n]) = (I, \omega)$ in $E^n(A)$. Let  $-\omega : (A/I)^n \to I/{I^2} $ denote the surjection given by $( \overline{f_1}, \dots, \overline{f_{n-1}}, -\overline{f_n} )$. Then we have the following proposition:

\begin{prop1} \label{prop1}
Let $A, P, I, \alpha, \omega$ and $-\omega$ be as above. Then $(I, \omega) + (I, -\omega) =0 $ in $E^n(A)$. 
\end{prop1}

\begin{proof}
We first show there is a unimodular row over $A$ which represents $(I, -\omega)$ in $E^n(A)$. Notice the $\sum_{0 \le i \le n } a_ip_i=0$ implies $\sum_{0 \le i \le n } a_if_i=0$. Let $\phi$ be the surjection: $A^{n+1} \to I $ given by $ ( f_0,  \dots, f_{n-1}, -f_n ) $ . Let $Q$ be the projective $A$-module of rank $n$ defined by 
$$Q=A^{n+1}/{(a_0e_0+ \dots + a_{n-1}e_{n-1}-a_ne_n)}.$$ 
Since $[a_0, a_1, \dots, a_n]$ is a unimodular row over $A$, there exist $b_i$'s $\in A$ such that $a_0b_0 + a_1b_1+ \dots + a_nb_n=1$. Let $q_n=e_n-(-b_n)(a_0e_0+ \dots +a_{n-1}e_{n-1} - a_ne_n)$ and $q_i=e_i-(b_i)(a_0e_0+ \dots +a_{n-1}e_{n-1} - a_ne_n)$ for $i= 0, \dots, n-1$, then $ Q= \sum_{0 \le i \le n } Aq_i $ , $\phi(q_n)=-f_n$ and $\phi(q_i)=f_i$ for $i= 0,  \dots, n-1 $. Thus the restriction of $\phi$ to $Q$ gives us a surjection: $Q \to I$, call it $\beta$. Then it is rather obvious to see that $e([a_0, \dots, a_{n-1}, -a_n]) = (I, -\omega) $ in $E^n(A)$ via the surjection $\beta$ and the projective $A$-module $Q$. 
\par
Next, we show the image of $[a_0, \dots, a_{n-1}, a_n] * [a_0, \dots, a_{n-1}, -a_n]$ under the group homomorphism $e$ is zero, where $*$ is the group operation on $$ Um_{n+1}(A)/{E_{n+1}(A)} $$ which is defined in [\ref {Vk2}] by Van der Kallen.  Applying [ \ref {Vk2}, Lemma 3.5 (i)], we have 
$$[a_0, \dots, a_{n-1}, a_n] * [a_0, \dots, a_{n-1}, -b_n] = 0 $$ 
in $ Um_{n+1}(A)/{E_{n+1}(A)} $. Also by [ \ref {Vk2}, Lemma 3.5 (v)], 
$$[a_0, \dots, a_{n-1}, -a_n]* [a_0, \dots, a_{n-1}, b_n^2]= [a_0, \dots, a_{n-1}, -a_nb_n^2].$$ 
But  $ [a_0, \dots, a_{n-1}, -a_nb_n^2] =[a_0, \dots, a_{n-1}, -b_n] $ in $Um_{n+1}(A)/{E_{n+1}(A)}$. Taking image in $E^n(A)$ under the group homomorphism $e$ and applying Corollary \ref {cor1}, it follows that  $ (I, \omega) + (I, -\omega) + 0 = 0 $ in $E^n(A)$. The proof of the proposition is complete.
\end{proof} 

From Lemma \ref {lm3} and Proposition \ref {prop1}, we have the following interesting corollary:

\begin{coro} \label{cor1'}
Let $A$ be a smooth affine $\mathbb{C}$-algebra of dimension $n+1 $($n \ge 4$, even). Let $[a_0, \dots,  a_n] \in Um_{n+1}(A)$. Then $e([a_0, \dots,  a_n])$ is 2-torsion in $ E^n(A)$. 
\end{coro}

\begin{proof}
Let $P$ denote the projective $A$-module defined by $[a_0, \dots,  a_n]$. Choose a general section of the dual of $P$, say $\alpha$, then $\alpha$ gives us a surjection: $P \to I$, where $I$ is a local complete intersection ideal of height $n$. From this surjection, we can write $e([a_0, \dots,  a_n]) = (I, \omega) $ in $E^n(A)$, where $\omega $ is some surjection: ${A/I}^n \to I/{I^2}$. By Lemma \ref {lm3}, we have $(I, \omega) = (I, -\omega)$ since -1 is a square in $A/I$. By Proposition \ref {prop1}, we have $2e([a_0, \dots,  a_n])=0$ in $E^n(A)$. The proof of the corollary is complete.
\end{proof}

\begin{lemma} \label {lm4}
Let $A$ be a commutative Noetherian ring containing a field $k$. Let $I$ be a proper ideal of height $n$ which is a local complete intersection in $A$, such that $I/{I^2}$ is free $A/I$-module of rank $n$.  Then there exists a regular sequence $f_1, \dots, f_n $ in $A$ and $s_1 \in I^2 $ such that
\begin{enumerate} 
\item $I=(f_1, \dots, f_n, s_1)$, $s_1(1-s_1) \in (f_1, \dots, f_n)$,  $I= (f_1, \dots, f_n) + I^2$, and
\item $ \{f_1, \dots, f_{n-1}, f_n - s_1^2 \}$  is a regular sequence in $A$.
\end{enumerate}
\end{lemma}

\begin{proof}
As in the proof of Lemma 2.3 in [\ref {Ma1}], we can find a regular sequence $f_1, \dots, f_n $ in $A$ such that $I=(f_1, \dots, f_n)+ I^2$. By [\ref {Mk}, Lemma 1], there exists $s \in I $ such that  $s(1-s) \in (f_1, \dots, f_n)$ and $I= (f_1, \dots, f_n, s)$. Since $s(1-s) \in (f_1, \dots, f_n)$, we may further assume that $s \in I^2 $. Notice that we can change $s$ by $\prod_{i=1}^{m} (s-b_if_n)$ for any positive integer $m$ and $b_i \in I$. If $p_1, \dots, p_t $ are the maximal elements in Ass$(A/{(f_1, \dots, f_{n-1})})$, then $f_n \notin p_1, \dots, p_t$. 

If $s \in p_1, \dots, p_t$, then $f_n -s^2 \notin p_1, \dots, p_t$, and we are through.

If, say for example, $s \notin p_1$, but $s+bf_n \in p_1$ for some $b \in I$, we replace $s$ by $s(s+bf_n)$ and assume $s \in p_1$. Repeating this procedure (that is, replacing $s$ by $\prod_{i=1}^{m} (s-b_if_n)$)  and reordering $p_i$ where $i \in \{1, \dots, t \}$ if necessary, we may assume that $s \in p_1, \dots, p_r $, $s-bf_n \notin p_k$ for $k >r$ and any $b \in I$.

Since $s \in p_1, \dots, p_r $, $f_n -s^2 \notin  p_1, \dots, p_r$. If $f_n -s^2 \notin  p_{r+1}, \dots, p_t$, then we are done. So by reordering $p_{r+1}, \dots, p_t$, we may assume $f_n -s^2 \notin p_{r+1}, \dots, p_{r+l}$ and $f_n -s^2 \in p_{r+l+1}, \dots, p_t$. Let $\lambda \in I \cap (\cap_{i=1}^{r+l} p_i ) \setminus \cup_{j=r+l+1}^t p_j $ (such $\lambda $ does exist), and $s_1= s+\lambda f_n$.  Then $f_n - s_1^2 = f_n - s^2 - \lambda f_n(2s + \lambda f_n)$, and $f_n - s_1^2 \notin p_1, \dots, p_{r+l}$ by our choice of $\lambda$. 

Now we claim that $f_n - s_1^2 \notin p_{r+l+1}, \dots, p_t$. If $f_n - s_1^2 \in p_j$ for some $j \in \{ r+l+1, \dots, t \}$, then $2s+\lambda f_n \in p_j$. Notice that since $A$ is a commutative ring containing a field $k$, either 2 is invertible in $A$ or 2 is zero in $A$. If 2 is zero in $A$, then $\lambda f_n \in p_j$, which is impossible.  If 2 is invertible in $A$, then $s+(1/2)\lambda f_n \in p_j$, which contradicts that $s-bf_n \notin p_k$ for $k >r$ and any $b \in I$. So the claim follows.

Therefor $f_n - s_1^2 $  is a nonzero divisor in $A/{(f_1, \dots, f_{n-1})}$. By our choice of $s_1$, we have that $I=(f_1, \dots, f_n, s_1)$, $s_1(1-s_1) \in (f_1, \dots, f_n)$, $s_1 \in I^2$, $I= (f_1, \dots, f_n) + I^2$, and $ \{f_1, \dots, f_{n-1}, f_n - s_1^2 \}$  is a regular sequence in $I$.
\end{proof}

The proof of the following lemma is a generalization of [\ref {zzm1}, Proposition 4.3] which is inspired by the statements of Bhatwadekar, Das, and Mandal in [\ref {Ma}, Lemma 6.1 and Proposition 6.2].

\begin{lemma} \label{lm5}
Let $A$ be a regular ring of dimension $d$ ($d \ge 3$) containing an infinite field $k$ and $n$ be an even integer such that $2n \ge d+3 $.. Let $I$ be an ideal of height $n$ such that 
$I=(f_1, \dots, f_{n-1}, f_n) + I^2$, where  $f_1, \dots, f_{n-1}, f_n$ form a regular sequence in $A$. Let 
$\omega_I$ be the surjection: $(A/I)^n \to I/I^2$ given by  $( \overline{f_1}, \dots, \overline{f_{n-1}}, \overline{f_n} )$ and $-\omega_I$ be the surjection: $(A/I)^n \to I/I^2$ given by  $( \overline{f_1}, \dots, \overline{f_{n-1}}, -\overline{f_n} )$.
Define
 $J=I^{(2)}=(f_1, \dots, f_{n-1}) + I^2 $.
 Then there exists a surjection $\omega$: $(A/{J})^n \to J/{J^2}$, such that
 $(I^{(2)}, \omega) = (I, \omega_I) + (I, -\omega_I) $
in $E^n(A)$.
\end{lemma}

\begin{proof}
Applying lemma \ref {lm4}, we can find $s \in A$ such that the image of $ f_n - s^2 $ in $A/{(f_1, \dots, f_{n-1})}$ is a nonzero divisor, $I=(f_1, \dots, f_n, s)$ and $s(1-s) \in (f_1, \dots, f_n)$. Let $K_1 = (f_1, \dots, f_n, 1-s)$, then $K_1 \cap I =(f_1, \dots, f_n)$. Since $\{f_1, \dots, f_{n-1}, f_n -s^2 \}$ is a regular sequence and $I= (f_1, \dots, f_{n-1}, f_n -s^2 ) + I^2$, we can write $(f_1, \dots, f_{n-1}, f_n -s^2 )=I \cap K_2$ for some $K_2$ in $A$, which is comaximal with $I$. If $K_1=A$, or $K_2=A$, then the conclusion of the lemma clearly holds.  So we may assume $K_1, K_2$ are ideals of  height $n$. Let $g=f_n - s^2$, then $gA + K_1 =A$, and hence $I, K_1, K_2$ are pairwise comaximal. It is clear that $g = -s^2$  is unit modulo $K_1$ and $f_n = s^2$ is unit modulo $K_2$. Since $I^{(2)}=(f_1, \dots, f_{n-1})+I^2$, $I^{(2)} \cap K_1\cap K_2 = (f_1, \dots, f_{n-1}, gf_n)$. So we have the surjective homomorphisms

 \[\begin {CD} A^n@>\alpha>> I\cap K_1@.\end{CD}  \]
\[\begin {CD} A^n@>\alpha'>> I\cap K_1@.\end{CD} \] 
 
given by  $( f_1, \dots, f_{n-1}, f_n )$ and $( f_1, \dots, f_{n-1}, -f_n )$ respectively. 
Then $\omega_I=\alpha \otimes A/I $,  $-\omega_I=\alpha' \otimes A/I $, 
Let $\omega_{K_1}=\alpha \otimes A/{K_1} $ and $-\omega_{K_1}=\alpha' \otimes A/{K_1} $.
Then we have 
$(I, -\omega_I) + (K_1 , -\omega_{K_1})=0 $ in $E^n(A)$.

Since $I^{(2)} \cap K_1\cap K_2 = (f_1, \dots, f_{n-1}, gf_n)$. So we also have two natural  surjective homomorphisms

 \[\begin {CD} A^n@>\beta>> I\cap K_2@.\end{CD}  \]
\[\begin {CD} A^n@>\gamma>> I^{(2)}\cap K_1 \cap K_2@.\end{CD} \]

given by $( f_1, \dots, f_{n-1}, g )$ and $( f_1, \dots, f_{n-1}, gf_n )$ respectively.

Let $\omega_{K_2}=\beta \otimes A/{K_2}$, and   $\omega=\gamma\otimes A/{I^{(2)}}$.
Since $g = -s^2$  is unit modulo $K_1$ and $f_n = s^2$ is unit modulo $K_2$,
from the surjections $\beta$ and $\gamma$ and lemma \ref {lm3}, we have:
$\beta \otimes A/I =\omega_I, \gamma \otimes A/{K_1}=-\bar s^2\omega_{K_1}, \gamma \otimes A/{K_2}= \bar s^2 \omega_{K_2}$ and the following relations in $E^n(A)$:
$(I, \omega_I) + (K_2 , \omega_{K_2})=0$,  and  $(I^{(2)}, \omega) + (K_1 , -\bar s^2\omega_{K_1})+ (K_2 , \bar s^2\omega_{K_2}) =(I^{(2)}, \omega) + (K_1 , -\omega_{K_1})+ (K_2 , \omega_{K_2})=0$. Hence $(I, \omega_I) + (I, -\omega_I) = (I, \omega_I) + (I, -\omega_I) + (I^{(2)}, \omega) + (K_1 , -\omega_{K_1})+ (K_2 , \omega_{K_2}) = (I, \omega_I)+(K_2 , \omega_{K_2})+(I^{(2)}, \omega)+ (I, -\omega_I)+(K_1 , -\omega_{K_1})= 0 + (I^{(2)}, \omega) + 0 = (I^{(2)}, \omega) $. Thus
 $(I^{(2)}, \omega) = (I, \omega_I) + (I, -\omega_I) $ in $E^n(A)$.
\end{proof}

Now we are ready to state our main theorem which gives an affirmative answer to Question \ref {qu}:

\begin{theorem} \label{thm1}
Let $A$ be a regular ring of dimension $d$ ($d \ge 3$) containing an infinite field $k$. Let $n$ be an even integer such that $2n \ge d+3 $. Let $I$ be an ideal in $A$ of height $n$ and $P$ be a projective $A$-module of rank $n$. Suppose $P \oplus A \approx A^{n+1}$ and there is a surjection $\alpha$: $P \to I$. Then $I$ is a set theoretic complete intersection ideal in $A$.
\end{theorem}

\begin{proof}
Let $\omega$ and $-\omega $ be as in Proposition \ref {prop1}. Then by Lemma \ref{lm5}, there exists a surjection $\omega'$: $(A/J)^n \to J/J^2$, where $J=I^{(2)}$, such that $(I^{(2)}, \omega') = (I, \omega_I) + (I, -\omega_I) $ in $E^n(A)$.  Using Proposition \ref {prop1}, we have $(I^{(2)}, \omega')=0 $ in $E^n(A)$. By Theorem 4.2 in [\ref {BS1}], $I^{(2)}$ is generated by $n$ elements and hence $I^{(2)}$ is a complete intersection ideal in $A$. Therefore, $I$ is a set theoretic complete intersection ideal in $A$.
\end{proof}

\begin{coro} \label{cor2}
Let $A$ be a smooth affine $\mathbb{C}$-algebra of dimension $n+1$, where $n \ge 4$ and even. Let $I$ be an ideal of height $n$. Suppose $I$ is the image of a stably free $A$-module of rank $n$. Then $I$ is a set theoretic complete intersection ideal in $A$.
\end{coro}

\begin{proof}
Applying cancellation theorem of Suslin [\ref{Su}, Theorem 1] and using Theorem \ref {thm1} ,  we see $I$ is a set theoretic complete intersection ideal in $A$. The proof of the corollary is complete.
\end{proof}

\begin{remark} \label{rmk1}
For the case when $n=3$ in Corollary \ref {cor2}, let $\alpha $ be the surjection: $ P \to I $, where $P$ is a stably free module of rank 3. By cancellation theorem of Suslin [\ref{Su}, Theorem 1], we have $P \oplus A \approx A^4$. By [\ref {Ba1}], $P$ has a unimodular element. Then using the same arguments as in Proposition \ref{prop0} except for applying cancellation theorem of Suslin instead of applying Bass cancellation, we can conclude that $I$ is a complete intersection ideal in $A$. Therefore, Corollary \ref {cor2} also holds when $n=3$. 
\end{remark}

The following corollary gives a positive answer to the Question \ref {qu0} in the case when $A$ is a smooth affine $\mathbb{C}$-algebra.

\begin{coro} \label {cor3}
Let $A$ be a smooth affine $\mathbb{C}$-algebra of dimension $n+1 $, where $n \ge 4$ and even. Let $I$ be a local complete intersection ideal of height $n$ such that $I/I^2$ is a free $A/I$-module of rank $n$. Suppose $(A/I)$ is torsion in $K_0(A)$. Then $I$ is a set theoretic complete intersection in $A$.
\end{coro}
\begin{proof}
Let $m$ be the integer such that $m(A/I)=0$ in $K_0(A)$. Write $I=(f_1, \dots, f_{n-1}, f_n) + I^2$, where $f_1, \dots, f_{n-1}, f_n$ form a regular sequence. Let  $ J=(f_1, \dots, f_{n-1}) + I^m $. By a result of Mandal [ \ref {Ma1}, Lemma 2.3], $(A/J)=m(A/I)$ in $K_0(A)$ and hence $(A/J)=0 $ in $K_0(A)$. By a result of Murthy [\ref {Mu1}, Theorem 2.2], we see $K$ is the image of a stably free $A$-module of rank $n$, where $ K=(f_1, \dots, f_{n-1}) + J^{(n-1)!} $. By Corollary \ref {cor2}, we conclude that $K$ is a set theoretic complete intersection in $A$ and hence so is $I$.
\end{proof}

\begin{remark} \label{rmk2}
As in Remark \ref {rmk1}, Corollary \ref {cor3} also holds when $n=3$
\end{remark}

\begin{acknowledgment}
I would like to take this opportunity to express my gratitude to my thesis adviser, Professor N. Mohan Kumar, for his guidance. Without his suggestions I would neither have thought of proving, nor would I have been able to prove these results in the field of commutative algebra. I also wish to thank Professor S. M. Bhatwadekar for many useful conversations about Euler Class group, without which this paper would not have taken its present form. 

Also, I would like to thank the referee for carefully going through the draft and for making very helpful suggestions. 
\end{acknowledgment}

\par

\bibliographystyle{amsplain}

\end{document}